\documentclass{amsart}
\usepackage{amsmath}

\newtheorem{thm}{Theorem}
\newtheorem{lem}[thm]{Lemma}

\theoremstyle{definition}
\newtheorem{defn}[thm]{Definition}

\newcommand{\CPb}{\overline{\mathbb{CP}}{}^{2}}
\newcommand{\CP}{{\mathbb{CP}}{}^{2}} 
\newcommand{\Z}{\mathbb{Z}}

\title{Exotic smooth structures on small 4-manifolds} 

\begin{document}

\author{Anar Akhmedov}
\address{School of Mathematics,
Georgia Institute of Technology,
Atlanta, GA, 30332-0160, USA}
\email{ahmadov@math.gatech.edu}

\author{B. Doug Park}
\address{Department of Pure Mathematics, 
University of Waterloo, 
Waterloo, ON, N2L 3G1, Canada}
\email{bdpark@math.uwaterloo.ca}

\dedicatory{Dedicated to Ronald J. Stern on the occasion of his sixtieth birthday}

\date{December 24, 2006.  Revised on February 23, 2007}


\begin{abstract} 
Let $M$\/ be either\/ $\CP\#3\CPb$ or\/ $3\CP\#5\CPb$. 
We construct the first example of a simply-connected irreducible 
symplectic $4$-manifold that is homeomorphic but not diffeomorphic to $M$.
\end{abstract}

\maketitle

\section{Introduction}

Given two 4-manifolds, $X$\/ and $Y$, we denote their connected sum by
$X\# Y$.  For a positive integer $m\geq 2$,  
the connected sum of $m$\/ copies of $X$\/ will be denoted by $mX$ for short. 
Let $\CP$ denote the complex projective plane and let $\CPb$ denote the underlying 
smooth $4$-manifold $\CP$ equipped with the opposite orientation.    
There is an extensive and growing literature on the 
construction of exotic smooth structures 
on the closed $4$-manifolds $\CP\#m\CPb$ and 
$3\CP\#n\CPb$ for some small positive integers $m$\/ and $n$.  
In the next two paragraphs, we briefly highlight some of the papers
that have appeared.  

The existence of an exotic smooth structure on a $4$-manifold was 
first proved by Donaldson in \cite{Do}.  Using SU(2) gauge theory, he 
showed that a Dolgachev surface $E(1)_{2,3}$ is 
homeomorphic but not diffeomorphic to $\CP\#9\CPb$.  Infinitely many
irreducible smooth structures on $\CP\#9\CPb$ were later constructed by
Fintushel, Friedman, Stern and Szab\'{o} in \cite{Fri,FS2,Sz2}. 
Kotschick showed in \cite{Ko1} 
that the Barlow surface in \cite{Ba} is
homeomorphic but not diffeomorphic to $\CP\#8\CPb$.  
More recently,
J. Park was able to construct an exotic $\CP\#7\CPb$ in \cite{Pa2} 
using the rational blowdown technique of Fintushel and Stern 
in \cite{FS1}.  Using a more general blowdown technique in \cite{Pa1},
Stipsicz and Szab\'{o} constructed an exotic $\CP\#6\CPb$ in \cite{SS2}.  
In \cite{FS3}, Fintushel and Stern constructed infinitely many irreducible 
smooth structures on $\CP\#n\CPb$ for $6\leq n\leq8$.  
Using Fintushel and Stern's technique in \cite{FS3}, 
J. Park, Stipsicz and Szab\'{o} were able to construct infinitely many
irreducible non-symplectic smooth structures on $\CP\#5\CPb$ in \cite{PSS}.  
The first \emph{symplectic}\/ 
exotic irreducible smooth structure on $\CP\#5\CPb$ was constructed in \cite{A2}
by the first author.  

The existence of infinitely many irreducible smooth structures on $3\CP\#19\CPb$
was shown by Friedman, Morgan, O'Grady, Stipsicz and Szab\'{o} 
in \cite{FM,MO,SS1}.   
Exotic smooth structures on $3\CP\#m\CPb$ for $14 \leq m\leq 18$ were first
constructed by Gompf in \cite{Go}.  The existence of infinitely many irreducible
smooth structures on $3\CP\#m\CPb$ for $14\leq m\leq 18$ were shown by 
Stipsicz, Szab\'{o} and Yu in \cite{St1,St2,Sz1,Yu}.  
Infinitely many irreducible smooth structures on
$3\CP\#m\CPb$ for $10\leq m\leq 13$ were constructed by the second author 
in \cite{P1,P2,P3}.  Using rational blowdown techniques, Stipsicz and Szab\'{o}
constructed infinitely many irreducible smooth structures on 
$3\CP\#9\CPb$ in \cite{SS3}, and J. Park constructed infinitely many irreducible smooth
structures on $3\CP\#8\CPb$ in \cite{Pa3}.   The first exotic irreducible 
smooth structure on
$3\CP\#7\CPb$ was constructed by the first author in \cite{A2}.   

In this paper, we use the techniques and constructions in \cite{A1} and \cite{A2} to
prove the following.   

\begin{thm}\label{thm:main}
Let $M$\/ be either\/ $\CP\#3\CPb$ or\/ $3\CP\#5\CPb$.  
There exists a smooth closed simply-connected irreducible symplectic\/ 
$4$-manifold that is homeomorphic but not diffeomorphic to\/ $M$.  
\end{thm}

The main technical difficulty in the proof of Theorem~\ref{thm:main} is the computation
of the fundamental groups of our exotic 4-manifolds.  Even though we cannot always 
completely pin down the fundamental groups at all stages of our construction, we are able to identify enough of the generators (up to conjugacy) and determine enough relations among them so that we are able to deduce, after careful choices of gluing, 
that the resulting 4-manifolds are simply-connected.  

In the forthcoming paper \cite{ABP} with R. \.{I}nan\c{c} Baykur, 
we will present some alternative constructions of an
irreducible symplectic smooth structure on $M$.  
In \cite{ABP}, we also construct infinitely many 
irreducible non-symplectic smooth structures on 
$M$, $3\CP\#7\CPb$ and other small $4$-manifolds.  

Here is how our paper is organized.  
Section~\ref{sec:def} contains some definitions and formulas that will be important throughout the paper.  Section~\ref{sec:blocks} 
quickly reviews the 4-manifolds that were constructed in
\cite{A1}.  These 4-manifolds will then serve as some of the building blocks in constructing an exotic $3\CP\#5\CPb$ in 
Section~\ref{sec:3CP-5CP} and an exotic $\CP\#3\CPb$
in Section~\ref{sec:CP-3CP}.
Finally, Section~\ref{sec:appendix} contains some technical 
computations which are crucial in showing that the fundamental groups
of our exotic 4-manifolds are trivial.  

\section{Generalized fiber sum}
\label{sec:def}

We first present a few standard definitions that will be used throughout the paper.   

\begin{defn} 
Let $X$\/ and $Y$\/ be closed oriented smooth\/ $4$-manifolds each containing a smoothly embedded surface $\Sigma$ of genus $g \geq 1$. Assume $\Sigma$ represents a homology of infinite order and has self-intersection zero in $X$\/ and $Y$, 
so that there exists a 
product tubular neighborhood, say $\nu\Sigma\cong \Sigma\times D^{2}$, in both $X$\/ and $Y$. Using an orientation-reversing and fiber-preserving diffeomorphism $ \psi : \Sigma\times S^1  \rightarrow \Sigma\times S^1$, we can glue $X \setminus\nu\Sigma$ and $Y\setminus \nu\Sigma$\/ along the boundary $\partial(\nu\Sigma)\cong 
\Sigma\times S^{1}$.  The resulting closed oriented smooth\/ $4$-manifold, denoted $X\#_{\psi}Y$, is called a \emph{generalized fiber sum}\/ of $X$\/ and $Y$\/ along $\Sigma$.  
\end{defn}

\begin{defn}
Let $e(X)$ and $\sigma(X)$ denote the Euler characteristic and the signature
of a closed oriented smooth\/ $4$-manifold $X$, respectively.  We define 
\begin{equation*}
c_1^2(X):= 2e(X)+3\sigma(X), \quad \chi_h(X):=\frac{e(X)+\sigma(X)}{4}.
\end{equation*}
\end{defn}

If $X$\/ is a complex surface, then $c_1^2(X)$ and $\chi_h(X)$ are 
the square of the first Chern class $c_1(X)$ and the holomorphic Euler 
characteristic, respectively.  Note that these can be used as the coordinates for the 
``geography problem'' for complex surfaces or irreducible smooth $4$-manifolds
(cf.\ \cite{GS}).  

For the simply-connected $4$-manifolds $m\CP\#n\CPb$, we have 
$e=2+m+n$\/ and $\sigma=m-n$.  Hence we get 
\begin{equation*}
c_1^2(m\CP\#n\CPb)=5m-n+4, \quad \chi_h(m\CP\#n\CPb)=\frac{m+1}{2}.
\end{equation*}

\begin{lem}\label{lemma:c_1^2 and chi} 
Let $X$ and $Y$ be closed, oriented, smooth\/ $4$-manifolds containing an embedded surface\/ $\Sigma$ of self-intersection\/ $0$. Then 
\begin{eqnarray*}
c_{1}^{2}(X\#_{\psi}Y) &=&  c_{1}^{2}(X) + c_{1}^{2}(Y) + 8(g-1),\\
\chi_{h}(X\#_{\psi}Y) &=&  \chi_{h}(X) + \chi_{h}(Y) + (g-1),
\end{eqnarray*}
where $g$ is the genus of the surface $\Sigma$. 
\end{lem}

\begin{proof}
The above formulas simply follow from the well-known formulas 
\begin{equation*}
e(X\#_{\psi}Y)= e(X) + e(Y) - 2e(\Sigma),\quad 
\sigma(X\#_{\psi}Y) = \sigma(X) + \sigma(Y).  \qedhere   
\end{equation*}
\end{proof}
 
If $X$\/ and $Y$\/ are symplectic $4$-manifolds and $\Sigma$ is a symplectic 
submanifold in both, then according to a theorem of Gompf (cf.\ \cite{Go}), 
$X\#_{\psi}Y$\/ admits a symplectic structure.  In such a case, we will call 
$X\#_{\psi}Y$\/ a \emph{symplectic sum}.

\section{Building blocks}
\label{sec:blocks}

We review the main construction in \cite{A1}.   
Let $K$\/ denote a trefoil knot in $S^3$.  
Let $\nu K$\/ denote the tubular neighborhood of $K$\/ in $S^3$.  
It is well-known (cf.\ \cite{BZ}) that 
\begin{equation*}
\pi_1(S^3\setminus \nu K) = \langle a, b \mid  aba = bab \rangle, 
\end{equation*}
and $b$\/ and $ab^2ab^{-4}$ 
represent a meridian and a longitude of $K$, respectively.  
It is also well-known that $\gamma_1=a^{-1}b$ and $\gamma_2=
b^{-1}aba^{-1}$ generate the image of the fundamental group of 
the genus one Seifert surface of $K$\/ under the inclusion-induced homomorphism.  
Let $M_K$ denote the result of 0-surgery on $K$.  Clearly, we have
\begin{eqnarray}
\pi_1(M_K\times S^1) &=& \pi_1(M_K) \oplus \Z \label{eq:pi_1(M_K x S^1)}\\
&=&  \langle a,b,x\mid aba=bab,\, ab^2ab^{-4}=1,\, [x,a]=[x,b]=1\rangle .
\nonumber
\end{eqnarray}
Since $K$\/ is a genus one fibered knot, 
$M_K$ is a $T^2$ fiber bundle over $S^1$ with a section $b$.
Thus $M_K\times S^1$ is a $T^2$ fiber bundle over $T^2$, with a section 
$S= b\times S^1=b\times x$.  From \cite{Th}, we know that there exists 
a symplectic form on $M_K\times S^1$ with respect to which 
both $S$\/ and a torus fiber $F$\/ are symplectic submanifolds.  
Let\/ ${\rm pr}:M_K\times S^1
\rightarrow S^1$ denote the projection map onto the second factor.  

\begin{lem}\label{lemma:C_S}
Let\/ $C_S=(M_K\times S^1)\setminus\nu S$ be the complement of 
a tubular neighborhood of a section\/ $S$ in $M_K\times S^1$.
Then we have 
\begin{equation*}
\pi_1(C_S) = \langle 
a, b, x \mid aba =bab, \, [x,a]=[x,b]=1
\rangle.  
\end{equation*}  
\end{lem}

\begin{proof}
Note that $C_S=(M_K\setminus \nu b)\times S^1=(S^3\setminus \nu K)
\times S^1$.  
\end{proof}

\begin{lem}\label{lemma:C_F}
Let\/ $C_F=(M_K\times S^1)\setminus\nu F$ be the complement of 
a tubular neighborhood of a fiber\/ $F$ in $M_K\times S^1$.  
Then we have 
\begin{eqnarray*}
\pi_1(C_F)&=&\langle \gamma_1',\gamma_2',d,y\mid
[\gamma_1',\gamma_2']=[y,\gamma_1']=[y,\gamma_2']=1,\\ 
&&d\gamma_1'd^{-1}=\gamma_1'\gamma_2',\,
d\gamma_2'd^{-1}=(\gamma_1')^{-1}\rangle.
\end{eqnarray*}
\end{lem}

\begin{proof}
$C_F$ is homotopy equivalent to a $T^2$ fiber bundle over a wedge of 
two circles.  The monodromy along the circle $y$\/ is trivial whereas
the monodromy along the circle $d$\/ is the same as the monodromy of $M_K$.   
Recall that $\gamma_1'=c^{-1}d$ and $\gamma_2'=d^{-1}cdc^{-1}$.
Thus we have 
\begin{eqnarray*}
d\gamma_1'd^{-1}&=&dc^{-1}dd^{-1}
= dc^{-1} = \gamma_1'\gamma_2',\\
d\gamma_2'd^{-1}&=&d d^{-1}cdc^{-1} d^{-1} 
= (cdc^{-1}) d^{-1}\\
&=& (d^{-1}cd) d^{-1} = d^{-1}c
= (\gamma_1')^{-1}.
\end{eqnarray*}
Note that $cdc=dcd$\/ implies $cdc^{-1}=d^{-1}cd$.
\end{proof}

\begin{lem}\label{lemma:pi_1(Y_K)}
Let $Y_K$ be the symplectic sum of two copies of $M_K\times S^1$,
identifying a section\/ $S$ in one copy with a fiber\/ $F$ in the other copy.  
If the gluing map
$\psi$ satisfies $\psi_{\ast}(x)=\gamma'_1$ and $\psi_{\ast}(b)=\gamma'_2$, then 
\begin{eqnarray*}
\pi_1(Y_K)&=&\langle 
a,b,x, \gamma_1',\gamma_2',d,y
\mid
aba=bab,\, [x,a]=[x,b]=1,\\
&& [\gamma_1',\gamma_2']=[y,\gamma_1']=[y,\gamma_2']=1,\, 
d\gamma_1'd^{-1}=\gamma_1'\gamma_2',\,
d\gamma_2'd^{-1}=(\gamma_1')^{-1},\\
&&x=\gamma_1',\, b=\gamma_2',\, ab^2ab^{-4} =[d,y]
\rangle \\
&=&\langle a,b,x,d,y
\mid aba=bab,\, [x,a]=[x,b]=1,\\
&&[y,x]=[y,b]=1,\, dxd^{-1}=xb,\, dbd^{-1}=x^{-1},\,
ab^2ab^{-4}=[d,y] \rangle.
\end{eqnarray*}
\end{lem}

\begin{proof}
By Seifert-Van Kampen Theorem, $\pi_1(Y_K)=(\pi_1(C_S)\ast\pi_1(C_F))/
\pi_1(T^3)$.  One circle factor of $T^3$ is identified with the longitude
of $K$ on one side and the meridian of the torus fiber in $M_K\times S^1$
on the other side.  This gives the last relation.  
\end{proof}

Inside $Y_K$, we can find a genus two symplectic submanifold 
$\Sigma_2$ which is the internal sum of a punctured fiber $F_0$ in $C_S$
and a punctured section $S_0$ in $C_F$.  The inclusion-induced homomorphism
maps the standard generators of $\pi_1(\Sigma_2)$ to $a^{-1}b$, 
$b^{-1}aba^{-1}$, $d$\/ and $y$\/ in $\pi_1(Y_K)$.  

\begin{lem}\label{lemma:Y_K complement}
There are nonnegative integers $m$ and $n$ such that
\begin{eqnarray}
\pi_1(Y_K\setminus \nu \Sigma_2 )&=&\langle
a,b,x,d,y;\, g_1,\dots,g_m \mid aba=bab,\, \label{eq:Y_K complement} \\
&&[y,x]=[y,b]=1,\, dxd^{-1}=xb,\, dbd^{-1}=x^{-1},\, \nonumber\\
&&ab^2ab^{-4}=[d,y], \,
r_1=\cdots=r_n=1,\, r_{n+1}=1 \rangle, \nonumber
\end{eqnarray}
where the generators $g_1,\dots,g_m$ and relators
$r_1,\dots,r_n$ all lie in the normal subgroup $N$\/ 
generated by the element\/ $[x,b]$, and the relator $r_{n+1}$
is a word in $x,a$ and elements of $N$.
Moreover, if we add an extra relation\/ 
$[x,b]=1$ to\/ $(\ref{eq:Y_K complement})$, then the relation\/ 
$r_{n+1}=1$ simplifies to\/ $[x,a]=1$.
\end{lem}

\begin{proof}
This follows from Seifert-Van Kampen Theorem.  Note that 
$[x,b]$ is a meridian of $\Sigma_2$ in $Y_K$.  Hence setting $[x,b]=1$
should turn $\pi_1(Y_K\setminus \nu \Sigma_2 )$ 
into $\pi_1(Y_K)$.   Also note that $a$\/ is another meridian of $K$, and so 
$a\times S^1$ is another torus section of the fiber bundle 
$M_K\times S^1$.
Thus $[x,a]$ is the boundary of a punctured torus section 
in $C_S\setminus \nu F_0$, and is no longer trivial in
$\pi_1(Y_K\setminus\nu\Sigma_2)$.  By setting
$[x,b]=1$, the relation $r_{n+1}=1$ is to turn into $[x,a]=1$.

It remains to check that the relations in $\pi_1(Y_K)$
other than $[x,a]=[x,b]=1$ remain the same in $\pi_1(Y_K\setminus\nu\Sigma_2)$.  
By choosing a suitable point $\theta\in S^1$ away from 
${\rm pr}|_{C_S}(F_0)$, 
the projection of 
the punctured fiber $F_0$ that forms a part of $\Sigma_2$, we obtain an embedding 
of the knot complement 
$(S^3\setminus\nu K)\times\{\theta\} \hookrightarrow 
C_S\setminus \nu F_0$.  This means that $aba=bab$\/ holds in 
$\pi_1(Y_K\setminus\nu\Sigma_2)$.   Since $[\Sigma_2]^2=0$, there 
exists a parallel copy of $\Sigma_2$ outside $\nu\Sigma_2$, wherein
the identity $ab^2ab^{-4}=[d,y]$ still holds.  The other four remaining relations in 
$\pi_1(Y_K)$ are coming from the monodromy of a torus fiber bundle over a
torus.  Since these four relations will now describe the monodromy of a 
punctured torus fiber bundle over a punctured torus, 
they still hold true in $\pi_1(Y_K\setminus\nu\Sigma_2)$.    
\end{proof}

Now take two copies of $Y_K$.  Suppose that the fundamental group of the second 
copy has $e,f,z,s,t$\/ as generators, and 
the inclusion-induced homomorphism in
the second copy maps the generators of $\pi_1(\Sigma_2)$ to $e^{-1}f$, 
$f^{-1}efe^{-1}$, $s$\/ and $t$\/ in $\pi_1(Y_K)$.  
Let $X_K$ denote the symplectic sum of two copies
of $Y_K$ along $\Sigma_2$, where the gluing map $\psi$ maps the generators as 
follows:
\begin{equation*}
\psi_{\ast}(a^{-1}b)=s,\; 
\psi_{\ast}(b^{-1}aba^{-1})= t,\; 
\psi_{\ast}(d)= e^{-1}f,\; 
\psi_{\ast}(y)=  f^{-1}efe^{-1}.
\end{equation*}

\begin{lem}\label{lemma:pi_1 X_K}
There are nonnegative integers $m$ and $n$ such that 
\begin{eqnarray*}
\pi_1(X_K)&=&\langle a,b,x,d,y;\;
e,f,z, s, t;\; g_1,\dots, g_m;\; h_1,\dots, h_m \mid \\
&&aba=bab,\, [y,x]=[y,b]=1,\\
&& dxd^{-1}=xb,\, dbd^{-1}=x^{-1},\,
ab^2ab^{-4}=[d,y],\\
&&r_1=\cdots=r_n=r_{n+1}=1,\, r_1'=\cdots=r_n'=r'_{n+1}=1 ,\, \\
&&efe=fef,\, [t,z]=[t,f]=1,\,\\
&& szs^{-1}=zf,\, sfs^{-1}=z^{-1},\,
ef^2ef^{-4}=[s,t],\\
&& d=e^{-1}f,\, y=f^{-1}efe^{-1},\, 
a^{-1}b=s,\, b^{-1}aba^{-1}=t,
[x,b]=[z,f] \rangle ,
\end{eqnarray*}
where $g_i,h_i$ $(i=1,\dots,m)$ and $r_j,r_j'$ $(j=1,\dots,n)$ all 
lie in the normal subgroup $N$ generated by
$[x,b]=[z,f]$.  Moreover, $r_{n+1}$ is a word in $x,a$ and elements of $N$, 
and $r'_{n+1}$ is a word in $z,e$ and elements of $N$.   
\end{lem}

\begin{proof}
This is just a straightforward application of Seifert-Van Kampen Theorem and
Lemma \ref{lemma:Y_K complement}.
\end{proof}

\begin{lem} \label{lemma:pi_1(X_K-Sigma_2)}
There are nonnegative integers $m,n,p$ and $q$ such that 
\begin{eqnarray*}
\pi_1(X_K\setminus \nu \Sigma_2 ) &=& 
\langle a,b,x,d,y;\,
e,f,z, s, t;\, g_1,\dots, g_m;\, h_1,\dots, h_m ;\, k_1,\dots, k_p \mid \\
&&aba=bab,\, [y,x]=[y,b]=1,\\
&& dxd^{-1}=xb,\, dbd^{-1}=x^{-1},\,
ab^2ab^{-4}=[d,y],\\
&&r_1=\cdots=r_n=r_{n+1}=1,\, r_1'=\cdots=r_n'=r'_{n+1}=1 ,\, \\
&&efe=fef,\, [t,z]=[t,f]=1,\,\\
&& szs^{-1}=zf,\, sfs^{-1}=z^{-1},\,
ef^2ef^{-4}=[s,t],\\
&& d=e^{-1}f,\, y=f^{-1}efe^{-1},\, 
a^{-1}b=s,\, b^{-1}aba^{-1}=t,\,\\
&&r''_1=\cdots=r''_q=1 \rangle ,
\end{eqnarray*}
where the elements $g_1,\dots,g_m,r_1,\dots,r_n$ lie in the normal subgroup 
$N$ generated by $[x,b]$, $h_1,\dots,h_m,r'_1,\dots,r'_n$ lie in the normal subgroup 
$N'$ generated by $[z,f]$, and $k_1,\dots,k_p,r''_1\dots,r''_q$ lie in the 
normal subgroup $N''$ generated by $[x,b][z,f]^{-1}$.  Moreover, 
$r_{n+1}$ is a word in $x,a$ and elements of $N$, 
and $r'_{n+1}$ is a word in $z,e$ and elements of $N'$.   
\end{lem}

\begin{proof}
Note that $[x,b][z,f]^{-1}$ is the meridian of $\Sigma_2$ in $X_K$.  
By Seifert-Van Kampen Theorem, setting $[x,b][z,f]^{-1}=1$ turns 
$\pi_1(X_K\setminus \nu\Sigma_2)$ into $\pi_1(X_K)$. 
\end{proof}

\section{Construction of an exotic $3\CP\#5\CPb$}
\label{sec:3CP-5CP}

In this section, we present a construction of a simply-connected, symplectic 
4-manifold $X$\/ that is homeomorphic but not diffemorphic to $3\CP\#5\CPb$. 
We will use the Seiberg-Witten invariant of $X$\/ to distinguish it from $3\CP\#5\CPb$.
Our first building block will be the symplectic $4$-manifold $X_{K}$ and a genus two symplectic submanifold $\Sigma_2\subset X_K$ from Section~\ref{sec:blocks}.  
The other building block will be the $4$-manifold $Y=T^4\#2\CPb$, the $4$-torus blown up twice.  $X$ will then be the 4-manifold obtained by taking the symplectic sum
of $X_{K}$ and $Y$\/ along $\Sigma_2$ and a certain genus two surface $\Sigma'_{2}\subset Y$. 

We find a symplectically embedded genus two surface $\Sigma'_{2}$ in $Y$\/ as follows.  
First we introduce the notation $T_{i,j}$ ($1\leq i < j\leq 4$)\/ for the 2-torus inside  $T^4$ that has nontrivial $i$th and $j$th circle factors.  For example, $T_{1,2}=S^1\times S^1\times \{{\rm pt}\}\times \{{\rm pt}\}$.  Let $p_{i,j}:T^4\rightarrow
T_{i,j}$ be the projection map.  Let $\omega_{i,j}$ be a standard product
volume form on $T_{i,j}$.  

Next we fix a factorization $T^4=T^{2}\times T^{2}$ and endow $T^4$ with a corresponding product symplectic form 
$\omega=p_{1,2}^{\ast}(\omega_{1,2})+p_{3,4}^{\ast}(\omega_{3,4})$.  
Consider one copy of a horizontal torus $T_{1,2}$, and one copy of a vertical torus $T_{3,4}$.  They are both symplectically embedded in $T^4$ with 
respect to $\omega$.  We symplectically resolve their intersection and obtain a symplectic surface of self-intersection $2$.  Next blow up twice to get a symplectic surface $\Sigma'_{2}$ of self-intersection $0$ in $Y=T^4\#2\CPb$. 

Note that the fundamental group of $Y$\/ is $\Z^4$. Let $\alpha_{i}$ ($i=1,\dots,4$)\/ denote the generators of $\pi_{1}(Y)$.  The fundamental group of the complement of a tubular neighborhood $\nu\Sigma'_{2}\cong\Sigma'_2 \times D^2$ is 
also $\Z^4$. It is also generated by the $\alpha_{i}$'s.  
This is because the normal circle $\mu'=\{{\rm pt}\}\times\partial D^2$ 
of $\nu\Sigma'_{2}$ (i.e. the meridian of 
$\Sigma'_2$) can be deformed into one of the exceptional spheres, and thus is 
trivial in $\pi_{1}(Y \setminus \nu\Sigma'_{2})$.   

Recall from Section~\ref{sec:blocks} 
(cf.\ \cite{A1}) that  $a^{-1}b$, $b^{-1}aba^{-1}$, $d$, $y$, and 
$\mu=\{{\rm pt}\}\times S^1=[x,b][z,f]^{-1}$ generate the 
inclusion-induced image of $\pi_1(\Sigma_{2}\times S^{1})$ 
inside $\pi_1(X_K\setminus\nu\Sigma_2)$.  
As before, let $\alpha_{1}$, $\alpha_{2}$, $\alpha_{3}$, $\alpha_{4}$ and $\mu'$ 
generate $\pi_1(\Sigma'_{2}\times S^{1})$.  
We choose the gluing diffeomorphism  $\psi : \Sigma_{2}\times S^1 \rightarrow 
\Sigma'_{2}\times S^1$ that maps the fundamental group generators as follows: 
\begin{equation*}
\psi_{\ast} (a^{-1}b) = \alpha_1,\;
\psi_{\ast} (b^{-1}aba^{-1}) = \alpha_2,\;
\psi_{\ast} (d) = \alpha_3 ,\; \psi_{\ast} (y) = \alpha_4 ,\;  \psi_{\ast} (\mu) = \mu' .
\end{equation*} 
By Gompf's theorem in \cite{Go}, 
$X := X_{K}\#_{\psi}(T^4\#2\CPb)$ is symplectic.  

\begin{lem}\label{lemma:pi_1(X)=1} 
$X$ is simply-connected. 
\end{lem}

\begin{proof}
By Seifert-Van Kampen Theorem, we have 
\begin{equation*}
\pi_{1}(X) \:=\: \frac{\pi_{1}(X_{K}\setminus \nu \Sigma_{2})\ast \pi_{1}(Y \setminus \nu \Sigma'_{2})}{\langle a^{-1}b = \alpha_{1},\, b^{-1}aba^{-1} = \alpha_{2},\, d = \alpha_{3},\, y = \alpha_{4},\, 
[x,b]=[z,f] \rangle}.
\end{equation*}
Note that the normal circle $\mu=[x,b][z,f]^{-1}$ of $X_{K} \setminus \nu\Sigma_{2}$ becomes trivial in $\pi_1(X)$ since $\mu'=1$.  This implies that the generators 
$k_1,\dots,k_p$ of $\pi_1(X_K\setminus\nu \Sigma_2)$ are trivial and the relations 
$r''_1=\cdots=r''_q=1$ are redundant in $\pi_1(X)$ 
(see Lemma~\ref{lemma:pi_1(X_K-Sigma_2)}).  
Since $\alpha_i$'s commute with one another, we get the following commutator relations in the fundamental group of $X$: 
$[a^{-1}b, b^{-1}aba^{-1}] = [a^{-1}b, d ] =  [a^{-1}b, y] =
[b^{-1}aba^{-1} , d] = [b^{-1}aba^{-1} , y] =
[d, y] =  1$.  
In summary we get the following presentation for the fundamental group of $X$.
\begin{eqnarray}
\pi_1(X)  &=& \langle 
a, b, x, d, y;\, e, f, z, s, t ;\, 
g_1,\dots, g_m;\, h_1,\dots, h_m \mid \label{eq:pi_1(X)}\\	
&& aba=bab,\, [y,x]=[y,b]=1,\nonumber\\
&& dxd^{-1}=xb,\, dbd^{-1}=x^{-1},\,
ab^2ab^{-4}=[d,y],\nonumber\\
&&r_1=\cdots=r_n=r_{n+1}=1,\, r_1'=\cdots=r_n'=r'_{n+1}=1 ,\,  \nonumber\\
&&efe=fef,\, [t,z]=[t,f]=1,\, \nonumber\\
&& szs^{-1}=zf,\, sfs^{-1}=z^{-1},\,
ef^2ef^{-4}=[s,t], \nonumber\\
&& d=e^{-1}f,\, y=f^{-1}efe^{-1},\, 
a^{-1}b=s,\, b^{-1}aba^{-1}=t,\, \nonumber\\
&& [s, t]=[s, d]=[s, y]=[t, d]=
[t, y]=[d, y] =1 \rangle. \nonumber
\end{eqnarray}
To prove $\pi_1(X)=1$, it is enough to prove that
$b=d=f=s=1$, since these will imply that all other generators
are trivial.  It is not hard to show that the following five identities hold in
$\pi_1(X)$.  The proof of these identities is postponed to the Appendix
(see Section~\ref{sec:appendix}).   
\begin{eqnarray}
dbd &=& bd^2b,\, \label{id:dbd}\\
bsb &=& sb^2s,\, \label{id:bsb}\\
sfs &=& fs^2f,\, \label{id:sfs}\\
fdf &=& df^2d,\, \label{id:fdf}\\
(bdb^{-1})s &=& s(bdb^{-1}). \label{id:[s,bdb^{-1}]=1}
\end{eqnarray}
Now rewrite (\ref{id:dbd}) as\/ $(d^{-1}b^{-1}d)(bdb^{-1}) = d$. Since $s$\/
commutes with $d$\/ and $bdb^{-1}$, we conclude that $s$\/ 
also commutes with $d^{-1}b^{-1}d$. It follows that $s$\/ commutes
with the inverse of $d^{-1}b^{-1}d$, i.e. 
$s(d^{-1}bd) = (d^{-1}bd)s$. 
Since $s$ commutes with $d$, we must have $sb = bs$.
Using $sb = bs$\/ in (\ref{id:bsb}), we get\/ $s = 1$. 
Then (\ref{id:sfs}) implies that\/ $f=1$.   Similarly, 
(\ref{id:fdf}) and (\ref{id:dbd}) in turn imply that $d = b = 1$. 
Thus $\pi_1(X)$ is trivial.
\end{proof}

From Lemma~\ref{lemma:pi_1(X_K-Sigma_2)}, 
it is easy to see that the abelianization of 
$\pi_1(X_K\setminus\nu\Sigma_2)$ is trivial.  By suitably adjoining an abelian 
group $\pi_1(Y\setminus\nu\Sigma'_2)\cong \Z^4$ to  
$\pi_1(X_K\setminus\nu\Sigma_2)$, 
we were able to infuse enough commutativity to 
eventually eliminate 
all the generators in the above proof.   Note that the relations 
$ab^2ab^{-4}=[d,y]$ and $ef^2ef^{-4}=[s,t]$ in (\ref{eq:pi_1(X)})
are not used in the above proof or in the proof of 
(\ref{id:dbd})--(\ref{id:[s,bdb^{-1}]=1}) 
in the Appendix.  An astute reader will notice that this and the fact that the
other meridians are eventually set to 1 allow us to
be a little lax on the orientation convention for meridians.  

\begin{lem} 
$e(X)=10$, $\sigma (X) = - 2$, $c_1^{2}(X) = 14$, and\/ $\chi_h(X) = 2$.
\end{lem}

\begin{proof} 
Let $Y=T^4\#2\CPb$ as before.  
We know that\/ $e(X)=e(X_K)+e(Y)+4$, 
$\sigma (X) = \sigma(X_{K}) + \sigma(Y)$, 
$c_1^{2}(X) = c_1^{2}(X_{K}) + c_1^{2}(Y) + 8$, and 
$\chi_{h}(X) = \chi_{h}(X_{K}) + \chi_{h}(Y) + 1$.  
We easily compute that $e(Y)=2$, 
$\sigma(Y)=-2$, 
 $c_1^{2}(Y)=-2$, 
and $\chi_{h}(Y)=0$.
Since $e(X_K)=4$, 
$\sigma(X_{K}) = 0$, 
$c_1^{2}(X_{K}) = 8$ 
and $\chi_{h} (X_{K}) = 1$ (cf.\ \cite{A1}), our results follow. 
\end{proof}

From Freedman's classification theorem (cf.\ \cite{Fre}) 
for simply-connected topological 4-manifolds, 
we conclude that $X$\/ is homeomorphic to $3\CP\#5\CPb$. 
It follows from Taubes's theorem (cf.\ \cite{Ta}) 
that ${\rm SW}_{X}( K_{X} ) = \pm 1$, where $K_X$ is the canonical class of $X$. 
Next we apply the connected sum theorem 
(cf.\ \cite{Wi}) for the Seiberg-Witten invariant to deduce that the Seiberg-Witten 
invariant is trivial for $3\CP\#5\CPb$. Since the Seiberg-Witten invariant is a diffeomorphism invariant, we conclude that $X$\/ is not diffeomorphic to 
$3\CP\#5\CPb$.  

Since $\pi_2(M_K\times S^1)=0$, we conclude that $M_K\times S^1$ is a
minimal symplectic 4-manifold.  Usher's theorem (cf.\ \cite{Us}) then implies that 
the symplectic sums $Y_K$ and $X_K$ are minimal as well.  
Both exceptional spheres $E_{1}$ and $E_{2}$ in $Y$\/ transversely intersect 
the genus two surface $\Sigma'_2$ once since $[\Sigma_2'] = 
[T_{1,2}] + [T_{3,4}] - [E_{1}] - [E_{2}]\in H_2(Y;\Z)$.  
Usher's theorem once again implies 
that the symplectic sum $X$\/ is a minimal symplectic 4-manifold.  
Since symplectic minimality implies irreducibility for simply-connected 
4-manifolds (cf.\ \cite{Ko2} for $b_2^+>1$ case), 
$X$\/ is also smoothly irreducible.

\section{Construction of an exotic $\CP\#3\CPb$}
\label{sec:CP-3CP}

In this section, we construct a simply-connected symplectic 4-manifold $U$\/ homeomorphic but not diffemorphic to $\CP\#3\CPb$. Using Usher's Theorem
(cf.\ \cite{Us}), we will distinguish $U$\/ from $\CP\#3\CPb$. 

The manifold $U$\/ will be the symplectic sum of the 4-manifold $Y_{K}$ in 
Section~\ref{sec:blocks} 
and $Q=(M_{K} \times S^{1})\#2\CPb$ along the genus two surfaces $\Sigma_{2}$ and 
$\Sigma''_{2}$.   The symplectic genus two submanifold $\Sigma_2''\subset Q$\/ 
is obtained by symplectically resolving 
the intersection of a torus fiber $F$\/ and a torus section $S$\/ of $M_K\times S^1$ 
(see Section~\ref{sec:blocks}) 
and then blowing up at two points.   

Let $g,h,z$\/ be generators of $\pi_1(M_K\times S^1)$ as in 
(\ref{eq:pi_1(M_K x S^1)}) such that the inclusion-induced homomorphisms
map the generators of $\pi_1(F)$ to $g^{-1}h$\/ and 
$h^{-1}ghg^{-1}$, and the generators of $\pi_1(S)$ to $h$\/ and $z$. 
Let us choose the gluing diffeomorphism  $\varphi : \Sigma_{2}\times S^1\rightarrow \Sigma''_{2}\times S^1$ that maps the elements $a^{-1}b$, $b^{-1}aba^{-1}$, $d$, 
$y$\/ and $\mu$\/ of $\pi_{1}(Y_{K} \setminus \nu \Sigma_{2})$ to the elements 
$g^{-1}h$, $h^{-1}ghg^{-1}$, $h$, $z$\/ and $\mu''$ of 
$\pi_{1}(Q \setminus \nu \Sigma_{2}'')$ according to the following rule: 
\begin{equation*}
\varphi_{\ast} (a^{-1}b) = h, \; \varphi_{\ast} (b^{-1}aba^{-1}) = z, \;
\varphi_{\ast} (d) = g^{-1}h, \; \varphi_{\ast} (y) = h^{-1}ghg^{-1}, \;
\varphi_{\ast} (\mu) = \mu'' .
\end{equation*} 
Here, $\mu$ and $\mu''$ denote the meridians of $\Sigma$ and $\Sigma_2''$.
It follows from Gompf's theorem (cf.\ \cite{Go}) that $U := Y_{K}\#_{\varphi}((M_{K} \times S^{1})\#2\CPb)$ is symplectic.  

\begin{lem}\label{lemma:pi_1(U)=1}
$U$\/ is simply-connected. 
\end{lem}

\begin{proof}
By Seifert-Van Kampen Theorem, we have 
\begin{equation*}
\pi_{1}(U) = \frac{\pi_{1}(Y_{K}\setminus \nu \Sigma_{2}) \ast \pi_{1}(Q \setminus \nu \Sigma''_{2})}{\langle a^{-1}b = h, \, b^{-1}aba^{-1} = z, \, d = g^{-1}h, \, y = h^{-1}ghg^{-1}, \, \mu = 1 \rangle}.
\end{equation*}
Since $\mu''$ is nullhomotopic in $Q \setminus \nu \Sigma_{2}''$, the normal circle $\mu = [x, b]$ of $\pi_{1}(Y_{K} \setminus \nu \Sigma_{2})$ becomes trivial in $\pi_1(U)$. This in turn implies that the generators $g_{1}, \dots , g_{m}$ of $\pi_{1}(Y_{K} \setminus \nu \Sigma_{2})$ become trivial in $\pi_1(U)$, the relations $r_1=\cdots=r_n=1$\/ become 
redundant, and the relation $r_{n+1}=1$ turns into $[x,a]=1$ 
(see Lemma~\ref{lemma:Y_K complement}). 
Since $\pi_{1}(Q\setminus \nu \Sigma''_{2})\cong\pi_1(Q)\cong\pi_1(M_K\times S^1)$,
we conclude that $\pi_1(U)$ is generated by  
\begin{equation*} 
a,\, b,\, x,\, d,\, y;\, g,\, h,\, z. 
\end{equation*}
The following relations hold in $\pi_1(U)$.  
\begin{eqnarray}
&&aba = bab, \, [x, b] = [x,a]=
 [y, x] = [y, b] = 1, \, \label{eq:pi_1(U) relations}\\
 && dxd^{-1} = xb, \, dbd^{-1} = x^{-1}, \, ab^{2}ab^{-4} = [d, y],\nonumber\\ 
 && ghg = hgh, \, gh^{2}gh^{-4} =1, \, [z, g] = [z, h] = 1, \nonumber\\
 && a^{-1}b = h, \, b^{-1}aba^{-1} = z, \, d = g^{-1}h, \, y = h^{-1}ghg^{-1} . \nonumber
\end{eqnarray} 
Using the above relations, we can obtain the following six identities in $\pi_1(U)$.  
The proof of these identities is postponed to the Appendix
(see Section~\ref{sec:appendix}).   
\begin{eqnarray}
bhb &=& hb^{2}h 		\label{id:bhb} \\
d(hbh^{-1}) &=& (hbh^{-1})d    \label{id:[d,hbh^{-1}]=1}\\
h(bdb^{-1}) &=& (bdb^{-1})h   \label{id:[h,bdb^{-1}]=1}\\
d^{2} &=& b^{-1}dbdb^{-1}     \label{id:d^2} \\
b^{2} &=& hbh^{-1}bh    \label{id:b^2} \\ 
h^{2} &=& dhd^{-1}hd	\label{id:h^2} 
\end{eqnarray}

Using (\ref{id:bhb}), we have $hbh^{-1} = b^{-1}hb^{2}$. Notice that by 
(\ref{id:[d,hbh^{-1}]=1}) 
$hbh^{-1}$ commutes with $d$. Then $d$ commutes with $b^{-1}hb^{2}$, so we have 
$d(b^{-1}hb^{2}) = (b^{-1}hb^{2})d$. 
This in turn implies that $(bdb^{-1})hb^{2} = hb^{2}d$. Since $h$\/ commutes with $bdb^{-1}$ by (\ref{id:[h,bdb^{-1}]=1}), we have $hbdb^{-1}b^{2} = hb^{2}d$. The last equality implies that $db = bd$. Using this fact in (\ref{id:d^2}), (\ref{id:b^2}) and (\ref{id:h^2}), 
we get $b = h = d = 1$. Now it easily follows from the relations of $\pi_{1}(U)$ that $a = x = y = g = z = 1$. Thus we have proved that $\pi_{1}(U)$ is trivial.   
\end{proof}

\begin{lem} 
$e(U)=6$, $\sigma (U) = - 2$, ${c_1}^{2}(U) = 6$, and $\chi_h(U) = 1$
\end{lem}

\begin{proof} 
Let $Q=(M_K\times S^1)\#2\CPb$.  
We have\/ $e(U)=e(Y_K)+e(Q)+4$, 
$\sigma (U) = \sigma(Y_{K}) + \sigma(Q)$, 
$c_1^{2}(U) = c_1^{2}(Y_{K}) + c_1^{2}(Q) + 8$, 
and $\chi_{h}(U) = \chi_{h}(Y_{K}) + \chi_{h}(Q) + 1$.  
Since $e(Q)=2$, $\sigma(Q)=-2$, $c_1^2(Q)=-2$, $\chi_h(Q)=0$, and $e(Y_K)= \sigma(Y_{K}) = c_1^{2}(Y_{K}) = \chi_{h} (Y_{K}) = 0$ (cf.\ \cite{A1}), 
our results follow.   
\end{proof}

By Freedman's theorem (cf.\ \cite{Fre}) and the lemmas above, we deduce that $U$\/ is homeomorphic to $\CP\#3\CPb$. Note that $U$\/ is a symplectic sum of a non-minimal 
4-manifold $Q = (M_{K}\times S^{1})\#2\CPb$ with a minimal 4-manifold $Y_{K}$.  Both exceptional spheres $E_{1}$ and $E_{2}$ in $Q$\/ transversely intersect 
the genus two surface $\Sigma''_2$ once since we have\/ 
$[\Sigma_2''] = [F] + [S] - [E_{1}] - [E_{2}]
\in H_2(Q;\Z)$.  It follows from Usher's Theorem (cf.\ \cite{Us}) 
that $U$\/ is a minimal symplectic 4-manifold.  
Since symplectic minimality implies irreducibility for simply-connected
4-manifolds (cf.\ \cite{HK} for $b_2^+=1$ case), 
$U$\/ is also smoothly irreducible.  
We conclude that it cannot be diffeomorphic to $\CP\#3\CPb$.

\section{Appendix}
\label{sec:appendix}

In this section we fill in the details left out in the proofs of Lemmas 
\ref{lemma:pi_1(X)=1} and \ref{lemma:pi_1(U)=1}
by proving the identities (\ref{id:dbd})--(\ref{id:[s,bdb^{-1}]=1}) in 
$\pi(X)$ and (\ref{id:bhb})--(\ref{id:h^2}) in $\pi_1(U)$.

\subsection*{Proof of (\ref{id:dbd})--(\ref{id:[s,bdb^{-1}]=1})}

From presentation (\ref{eq:pi_1(X)}), we get $x=db^{-1}d^{-1}$.  Substituting this into
the relation $dxd^{-1}=xb$, we get $d^2b^{-1}d^{-2}=db^{-1}d^{-1}b$, which can be 
rearranged to $b^{-1}d^{-2}b^{-1}=d^{-1}b^{-1}d^{-1}$.  By taking inverses of both 
sides, we get (\ref{id:dbd}).  
Similarly, substituting $z=sf^{-1}s^{-1}$ into the relation $szs^{-1}=zf$\/ yields
$s^2f^{-1}s^{-2}=sf^{-1}s^{-1}f$.  Rearrange this into 
$f^{-1}s^{-2}f^{-1}=s^{-1}f^{-1}s^{-1}$ and take inverses to obtain 
(\ref{id:sfs}).

Next note that $a=bs^{-1}$, and thus $aba=bab$\/ implies that
$bs^{-1}b^{2}s^{-1}=b^2 s^{-1}b$.  Hence 
$b^{-1}s^{-1}b^{2}s^{-1}b^{-1}=s^{-1}$, and by taking the inverses of both sides,
we obtain 
\begin{equation}\label{id:s}
(bsb^{-1})(b^{-1}sb) =s.
\end{equation}
Now $t=b^{-1}aba^{-1}=s^{-1}bsb^{-1}$.  
Since $st=ts$, we have $bsb^{-1}=s^{-1}bsb^{-1}s$.
From (\ref{id:s}), we get $b^2sb^{-2}s=bsb^{-1}$.
Combining the last two identities, we obtain
$b^2sb^{-2}s=s^{-1}bsb^{-1}s$,  which simplifies to $b^2sb^{-1}=s^{-1}bs$.
The last identity is easily changed to (\ref{id:bsb}).

Similarly, note that $e=fd^{-1}$, and thus $efe=fef$\/ implies that
$fd^{-1}f^2 d^{-1}=f^2 d^{-1}f$.  Then $f^{-1}d^{-1}f^2d^{-1}f^{-1}=d^{-1}$, and so
\begin{equation}\label{id:d}
fdf^{-2}df=d.
\end{equation}
Now $y=f^{-1}efe^{-1}=d^{-1}fdf^{-1}$.  Since $dy=yd$, we have 
$fdf^{-1}=d^{-1}fdf^{-1}d$.
From (\ref{id:d}), we get $f^2df^{-2}d=fdf^{-1}$.  
Combining the last two identities, we obtain
$f^2df^{-2}d=d^{-1}fdf^{-1}d$, which simplifies to $f^2df^{-1}=d^{-1}fd$.  
The last identity is easily changed to (\ref{id:fdf}).

Finally, since $dt=td$, we have $ds^{-1}bsb^{-1}=s^{-1}bsb^{-1}d$.
Since $d$\/ commutes with $s$\/ and $s^{-1}$, we conclude that 
\begin{equation*}
d(bsb^{-1})=(bsb^{-1})d.
\end{equation*}
From (\ref{id:s}), we also know that $d$\/ commutes with the product
$(bsb^{-1})(b^{-1}sb)$.  It follows that $d$\/ commutes with $b^{-1}sb$, i.e.
\begin{equation}\label{id:[d,b^{-1}sb]=1}
d(b^{-1}sb)=(b^{-1}sb)d.  
\end{equation}
Now it is easy to see that (\ref{id:[d,b^{-1}sb]=1}) can be rearranged
to (\ref{id:[s,bdb^{-1}]=1}).

\subsection*{Proof of (\ref{id:bhb})--(\ref{id:h^2})}

From (\ref{eq:pi_1(U) relations}), $x=db^{-1}d^{-1}$.  Plugging this into $dxd^{-1}=xb$, 
we get $d^2b^{-1}d^{-2}=db^{-1}d^{-1}b$.  This implies that $d^{-2}=
bd^{-1}b^{-1}d^{-1}b$, which gives (\ref{id:d^2}).  

Next we have $a^{-1}b=h$, and hence $a=bh^{-1}$.
Plugging this last identity into $aba=bab$, we get $bh^{-1}b^2h^{-1}=
b^2h^{-1}b$.  This simplifies to $h^{-1}b^2h^{-1}=b h^{-1} b$, which can be
easily rearranged to (\ref{id:b^2}).

Similarly, we have $d=g^{-1}h$, and hence $g=hd^{-1}$.  Plugging the last identity into $ghg=hgh$, we get $hd^{-1}h^2d^{-1}=h^2d^{-1}h$.  
It follows that $d^{-1}h^2d^{-1}=hd^{-1}h$, which can be easily arranged to 
(\ref{id:h^2}).

The relation 
$h = a^{-1}b$\/ gives $a=bh^{-1}$.  Thus $z = b^{-1}aba^{-1}=h^{-1}bhb^{-1}$.  
Since $[z, h] = 1$, we get 
$h^{-1}bhb^{-1}h = bhb^{-1}$. This simplifies to 
\begin{equation}\label{id:[h,bhb^{-1}]=1}
(bhb^{-1})h = h(bhb^{-1}).
\end{equation} 
From (\ref{id:b^2}), we have $b^2 = hbh^{-1}bh$, which can be rearranged as 
$hb^{-1}h^{-1}b = bhb^{-1}$.  Using (\ref{id:[h,bhb^{-1}]=1}), 
we have $bhb^{-1} = h^{-1}(bhb^{-1})h = b^{-1}h^{-1}bh$. 
The identity $bhb^{-1}=b^{-1}h^{-1}bh$\/ can be easily rearranged to
(\ref{id:bhb}).

The relation $g^{-1}h = d$\/ gives $g=hd^{-1}$.  Thus $y=
h^{-1}ghg^{-1}=d^{-1}hdh^{-1}$.  Since 
$[y, b] = 1$, we conclude that  
\begin{equation}\label{id:d^{-1}(hdh^{-1})b}
d^{-1}(hdh^{-1})b = bd^{-1}(hdh^{-1}) . 
\end{equation}
Now (\ref{id:h^2}) implies that $hdh^{-1} =
dh^{-1}d^{-1}h$. 
Substituting this last identity into (\ref{id:d^{-1}(hdh^{-1})b}), 
we obtain $h^{-1}d^{-1}hb=bh^{-1}d^{-1}h$.  This implies that 
$d^{-1}(hbh^{-1})=(hbh^{-1})d^{-1}$, which can be rearranged to 
(\ref{id:[d,hbh^{-1}]=1}).

Finally, note that $z$\/ commutes with both $g$\/ and $h$, so it must commute 
with $d = g^{-1}h$.  
In the proof of (\ref{id:bhb}) we saw that $z = h^{-1}bhb^{-1}$.  Substitute this 
into $dz = zd$ and obtain 
\begin{equation}\label{id:d(h^{-1}bhb^{-1})}
d(h^{-1}bhb^{-1}) = (h^{-1}bhb^{-1})d .
\end{equation}
From (\ref{id:b^2}), we have $h^{-1}bh=b^{-1}h^{-1}b^2$. 
Substituting this into 
(\ref{id:d(h^{-1}bhb^{-1})}), we obtain $db^{-1}h^{-1}b = b^{-1}h^{-1}bd$. 
This can be rearranged to $(bdb^{-1})h^{-1} =
h^{-1}(bdb^{-1})$, which implies (\ref{id:[h,bdb^{-1}]=1}).

\section*{Acknowledgments}
The authors thank R. \.{I}nan\c{c} Baykur
and Ronald J. Stern for helpful discussions.  
Some of the computations have been double-checked by 
using the computer software package GAP (cf.\ \cite{GAP}).  
A. Akhmedov was partially supported by the NSF grant FRG-0244663.  
B. D. Park was partially supported by CFI, NSERC and OIT grants.

\end{document}